# Instant Evaluation and Demystification of $\zeta(n), L(n,\chi)$ that Euler, Ramanujan Missed – I

THE INSTITUTE OF SCIENCE,
15, MADAM CAMA ROAD,
MUMBAI-400 032,
INDIA.

e-mail address : v_v_rane@yahoo.co.in

**Abstract** : For the Hurwitz zeta function $\zeta(s,\alpha)$, as a function of second variable $\alpha$ we obtain its Taylor series expansions about various points and show that they are all polynomials unlike $\frac{\partial}{\partial s}\zeta(s,1+\alpha)$ when s is a non-positive integer. On using functional equations, this results in explicit evaluation of $\zeta(2n)$ and of $L(n,\chi)$ when integral $n \geq 1$ and $\chi$ are both odd or both even. Here $\zeta(s)$ is Riemann zeta function and $L(s,\chi)$ is Dirichlet L-series. This also explains why $\zeta(2n+1)$ or $L(n,\chi)$ with $n \geq 1$ and $\chi$, of opposite parity cannot be evaluated explicitly. We also generalise Dirichlet's expression for $L(1,\chi)$ with real Dirichlet character, to any general Dirichlet character. We also deal with Lerch's zeta function on similar lines.

# Instant Evaluation and Demystification of $\zeta(n), L(n,\chi)$ that Euler, Ramanujan Missed - I

V.V.RANE
THE INSTITUTE OF SCIENCE,
15, MADAM CAMA ROAD,
MUMBAI-400 032,
INDIA.

e-mail address : v_ v_ rane @yahoo.co.in

Let $s = \sigma + it$ be a complex variable, where $\sigma$ and $t$ are real. For $\alpha \neq 0, -1, -2, -3, \ldots$, let $\zeta(s, \alpha)$ be the Hurwitz zeta function defined by $\zeta(s,\alpha) = \sum_{n \geq 0}(n+\alpha)^{-s}$ for $\sigma > 1$ and its analytic continuation. Thus $\zeta(s,1) = \zeta(s)$, the Riemann zeta function. For any integer $k \geq 1$, let $\zeta_k(s,\alpha) = \sum_{n \geq k}(n+\alpha)^{-s}$ for $\sigma > 1$ and its analytic continuation so that $\zeta_k(s,\alpha) = \zeta(s, k+\alpha) = \zeta(s,\alpha) - \sum_{0 \leq n \leq k-1}(n+\alpha)^{-s}$. For an integer $q \geq 1$, let $\chi(\bmod q)$ be a Dirichlet character and let $L(s,\chi)$ be the corresponding Dirichlet L-series. Thus $L(s,\chi) = q^{-s} \cdot \sum_{a=1}^{q}\chi(a)\zeta(s, \tfrac{a}{q})$. Let $\tau(\chi) = \sum_{a=1}^{q}\chi(a)e^{\frac{2\pi i a}{q}}$. For real $\lambda$ and for complex $\alpha \neq 0, -1, -2, \ldots$, let $\phi(\lambda, \alpha, s)$ be Lerch's zeta function defined by

$\phi(\lambda, \alpha, s) = \sum_{n \geq 0} e^{2\pi i \lambda n}(n+\alpha)^{-s}$ for $\sigma > 1$ and its analytic continuation. Note that if $\lambda$ is an integer then $\phi(\lambda, \alpha, s) = \zeta(s, \alpha)$. If $\lambda$ is not an integer, it is well-known that Lerch's zeta function $\phi(\lambda, \alpha, s)$ is an entire function of the complex variable s. We write $\phi(\lambda, s) = \sum_{n \geq 1} e^{2\pi i \lambda n} n^{-s}$ for $\sigma > 1$ and its analytic continuation. Note that $\phi(\lambda, 1, s) = e^{-2\pi i \lambda}\phi(\lambda, s)$. Thus $\phi(\lambda, s)$ is an entire function of s, if $\lambda$ is not an integer. Also note that $\frac{\partial^n}{\partial \alpha^n}\phi(\lambda, \alpha, s) = (-1)^n s(s+1)\ldots\ldots(s+n-1)\phi(\lambda, \alpha, s+n)$ for any integer $n \geq 0$. In what follows, m,n shall denote non-negative integers. Also $\zeta'(s,\alpha), L'(s,\chi)$ shall denote derivatives with respect to s, of the functions concerned. We write $S(s,\chi) = \sum_{n=1}^{q}\chi(n)n^{-s}$. In particular, we have for any integer $m \geq 0$, $S(-m, \chi) = \sum_{a=1}^{q}\chi(a)a^m$. In what follows, $\Gamma$ shall denote the gamma function.

It is well-known that $\zeta(2n) = \frac{(-1)^{n+1}2^{2n}B_{2n}\pi^{2n}}{2(2n)!}$ for any integer $n \geq 1$. Here, $B_n$'s are the pre-defined Bernoulli numbers, which are rational numbers. That is, $\zeta(2n)$ is a rational multiple of $\pi^{2n}$. However, $\zeta(2n+1)$ for $n \geq 1$ is not computable explicitly. On the other hand, if $L(s) = \sum_{n \geq 0} \frac{(-1)^n}{(2n+1)^s}$ for $\sigma > 1$ and its analytic continuation, then it known that $L(2n+1)$ for $n \geq 0$ is explicitly computable and $L(2n)$ for $n \geq 1$ cannot be computed explicitly. Note that $L(s) = \sum_{n \geq 1} \chi_4(n) n^{-s}$ for $\sigma > 1$, where $\chi_4$ is the unique odd character mod 4. Similarly, if $L_1(s) = \sum_{n \geq 1} \chi_3(n) n^{-s}$ for $\sigma > 1$ and its analytic continuation, where $\chi_3$ is the unique odd character mod 3, $L_1(2n+1)$ for $n \geq o$ is explicitly computable and $L_1(2n)$ for $n \geq 1$ cannot be computed explicitly. It is also known that $L(n, \chi)$ can be evaluated in terms of generalised Bernoulli polynomials, if $n \geq 1$ and $\chi(\mod q)$ are both even or both odd. However, $L(n, \chi)$ cannot be explicitly evaluated, if $n > 1$ and $\chi(\mod q)$ are of opposite parity. However, $L(1, \chi)$ can be evaluated for entirely different reasons. The question is as to why such things happen. In effect, it means that if we were not aware of Bernoulli polynomials and the explicit expressions thereof or if we were not aware of Bernoulli numbers and their values, we would not have been able to evaluate $\zeta(2n)$ for $n \geq 1$ or $L(n, \chi)$, when $n > 1$ and $\chi(\mod q)$ are of the same parity. The important question is whether evaluation of $\zeta(n)$ and $L(n, \chi)$ can be delinked from the theory of Bernoulli numbers, Bernoulli polynomials and generalized Bernoulli numbers and polynomials.

In author [3], it has been shown that for a fixed complex number s, $\zeta_1(s, \alpha) - \zeta(s)$ is an analytic function of $\alpha$ in the unit disc $|\alpha| < 1$ of complex plane with its power series expression $\zeta(s, 1+\alpha) - \zeta(s) = \sum_{n \geq 1} \frac{(-\alpha)^n s(s+1).....(s+n-1) \zeta(s+n)}{n!}$. The most important fact that emerges in author [3] (which we can see in what follows) is that for a fixed complex number s, $\zeta(s, \alpha) - \zeta(s)$ is an analytic function of $\alpha$ in the whole complex plane except for $\alpha = o, -1, -2, ....................$ However, it later turned out that Ramanujan [2] had already stated (without proof) the power series expansion of $\zeta(s,1-\alpha)$ for real $\alpha$ in the form $\zeta(s,1-\alpha) = \zeta(s) + \sum_{n\geq 1} \frac{\alpha^n}{n!} s(s+1)\ldots(s+n-1)\zeta(s+n)$ in a casual fashion. Surprisingly, Ramanujan does not elaborate on this power series expression nor does he give its origin, nor does he give the power series expression of $\zeta(-n,\alpha)$ for any integer $n \geq 0$. Actually, it can be shown that the power series of $\zeta(-n,\alpha)$ is a polynomial in $\alpha$. On the other hand, Ramanujan gives the power series expression of $\zeta'(-n,1+\alpha)$ for any integer $n \geq 0$. See entry 28 of Berndt [1]. It seems, the expression for $\zeta'(-n,1+\alpha)$ in entry 28 of Berndt [1] needs correction. For this see author [4]. Though Ramanujan gives the power series expression of $\zeta'(-n,1+\alpha)$, he does not give the power series expression of $\zeta'(s,1+\alpha)$ for any general complex number s nor does he deal with any possible connection between the power series expressions of $\zeta(s,1+\alpha)$ and $\zeta'(s,1+\alpha)$ nor does he use his power series expression for $\zeta'(-n,1+\alpha)$ to give rapidly convergent series expression for $\zeta(2n+1)$ for $n \geq 1$; or for $L(n,\chi)$, when $n \geq 1$ and $\chi(\mod q)$ are of opposite parity. Actually, it can be shown that the power series of $\zeta'(s,1+\alpha)$ can be obtained from the power series of $\zeta(s,1+\alpha)$ by term-by-term differentiation with respect to s. For this see author [4].

Exactly on the same lines as $\zeta(s,\alpha)$, if $\lambda$ is real, we get $\phi(\lambda,\alpha,s)$ as an analytic function of the complex variable $\alpha$ for $0 < |\alpha| < 1$ for a fixed complex number s and have the power series $\phi(\lambda,\alpha,s) - \alpha^{-s} = \phi(\lambda,s) + \sum_{n\geq 1} \frac{(-\alpha)^n}{n!} s(s+1)\ldots(s+n-1)\phi(\lambda,s+n)$. This follows by imtating the proof of author [3] in the case of $\zeta(s,\alpha)$ and on noting that $|\phi(\lambda,\alpha,s)| \leq \zeta(\sigma)$ for $\sigma > 1$ and $\lambda$ real. Also we can see that $\phi(\lambda,\alpha,s)$ is an analytic function of the complex variable $\alpha$ except for $\alpha = 0,-1,-2,\ldots$, exactly on the same lines as in the case of $\zeta(s,\alpha)$. We can also give power series expression in $\alpha$ corresponding to $\phi(\lambda,\alpha,s)$ and evaluate $\phi(\lambda,\alpha,-m)$, if m is a non-negative integer. See our Theorem 4 below.

Our approach can also be used to study $Z_r(s,\alpha)$, where for integral $r \geq 1$, and for $0 < \alpha \leq 1$, $Z_r(s,\alpha) = \sum_{n_1,n_2,\ldots,n_r \geq 0}(n_1 + n_2 + \ldots + n_r + \alpha)^{-s}$ for $\sigma > r$ and its analytic continuation. Note that $Z_1(s,\alpha) = \zeta(s,\alpha)$. Also note that $Z_r(s,\alpha) = \sum_{n \geq 0}\binom{n+r-1}{r-1}(n+\alpha)^{-s}$. We write $Z_r(s,1) = \zeta_r(s)$. Note $\zeta_1(s) = \zeta(s)$; $\zeta_2(s) = \zeta(s-1)$; and for $r \geq 3$, $\zeta_r(s) = \frac{1}{(r-1)!}\sum_{j=1}^{r-1}A_j\zeta(s-j)$, where $A_j = A_j(r)$'s are elementary symmetric sums of the numbers $1,2,3,\ldots,r-2$. As in the case of $\zeta(s,\alpha)$ and $\phi(\lambda,\alpha,s)$ we can give a power series expression of $Z_r(s,\alpha) - \alpha^{-s}$ for $|\alpha| < 1$ and also give the polynomial expression of $Z_r(-m,\alpha)$ for any integer $m \geq 0$. See our Theorem 5 below.

Our approach can also be used to study multiple Hurwitz zeta values $\zeta(s_1,s_2,\ldots,s_r;\alpha)$, where for a fixed integer $r \geq 1$ and for complex variables $s_1,s_2,\ldots,s_r$, we have

$\zeta(s_1,s_2,\ldots,s_r;\alpha) = \sum_{n_1 \geq 0}(n_1+\alpha)^{-s_1}\sum_{n_2 > n_1}(n_2+\alpha)^{-s}\ldots\sum_{n_r > n_{r-1}}(n_r+\alpha)^{-s_r}$ in some right half planes of complex variables $s_1,s_2,\ldots,s_r$ and its analytic continuation. Note that

$\frac{\partial}{\partial\alpha}\zeta(s_1,s_2,\ldots,s_r;\alpha) = -s_1\zeta(s_1+1,s_2,\ldots,s_r;\alpha) - s_2\zeta(s_1,s_2+1,\ldots,s_r;\alpha) - \ldots - s_r\zeta(s_1,s_2,\ldots,s_r+1;\alpha)$.

It is to be stressed that what we have given here is 'power series approach'. We can obtain similar results by using 'Integration approach' also.

A simple consequence of the power series expression for $\zeta_1(s,\alpha)$ is as follows, which we state as Proposition 1.

**Proposition 1** : We have for any integer $k \geq 1$ and for any fixed complex number s and for any complex number $\alpha$ with $|\alpha| < k$,

$$\zeta(s,\alpha) = \sum_{0 \leq n \leq k-1}(n+\alpha)^{-s} + \left(\zeta(s) - \sum_{1 \leq \ell \leq k-1}\ell^{-s}\right) + \sum_{n \geq 1}\frac{(-\alpha)^n}{n!}s(s+1)....(s+n-1)\left(\zeta(s+n) - \sum_{1 \leq \ell \leq k-1}\ell^{-s-n}\right)$$

**Remark**: The proof of Proposition 1, follows on considering $\zeta_k(s,\alpha) = \sum_{n \geq k}(n+\alpha)^{-s}$ for $\sigma > 1$ and then proceeding exactly on the same lines as in author [3], where we have dealt with the case $k=1$.

**Corollary of Proposition 1:** For any non-principal Dirichlet character $\chi(\bmod q)$ and for integer $k \geq 1$, we have

$$L(s,\chi) = \sum_{1 \leq n \leq q(k-1)}\chi(n)n^{-s} + \sum_{n \geq 1}\frac{(-1)^n}{n!}q^{-s-n}s(s+1).....(s+n-1)\left(\zeta(s+n) - \sum_{o \leq \ell \leq k-1}\ell^{-s-n}\right)\left(\sum_{a=1}^{q}\chi(a)a^n\right)$$

**Note**: This follows from Proposition 1 in view of $L(s,\chi) = q^{-s}\sum_{a=1}^{q}\chi(a)\zeta(s,\frac{a}{q})$.

**Remark**: From the Proposition 1, as $k \geq 1$ is an arbitrary integer, it is clear that $\zeta(s,\alpha)$ is a continuous function of $\alpha$ in the whole complex plane, if s is fixed with $\text{Re}\, s < o$; and is an entire function of $\alpha$, if s is a fixed non-positive integer and is an analytic function of $\alpha$ in the whole complex plane except for $\alpha = o, -1, -2, ..............$, if s is any fixed complex number.

One of our main objects is to show that the power series of $\zeta_1(-m,\alpha) = \zeta(-m, 1+\alpha)$ for any integer $m \geq 0$ is a polynomial of degree m+1 with rational coefficients and consequently the power series of $\zeta(-m,\alpha) = \alpha^m + \zeta(-m, 1+\alpha)$ is a polynomial of degree m+1 with rational coefficients. In particular, $\zeta(-m,\alpha)$ is rational, if $\alpha$ is rational. In fact, we prove the following Theorem 1.

**Theorem 1**: We have for any integer $m \geq 0$ and for any complex number $\alpha$,

$$\zeta(-m,\alpha) = \sum_{k=0}^{m}\binom{m}{k}\zeta(-k)\alpha^{m-k} + \alpha^m - \frac{\alpha^{m+1}}{m+1}, \text{ where } \binom{0}{0} = 1.$$

**Corollary 1**: If $k \geq 1$ is any integer, then $\zeta(-m,k) = \sum_{\ell=0}^{m}\binom{m}{\ell}\zeta(-\ell)k^{m-\ell} + k^m - \frac{k^{m+1}}{m+1}$.

Noting that $L(s,\chi) = q^{-s}\sum_{a=1}^{q}\chi(a)\zeta(s,\frac{a}{q})$ so that for any integer $m \geq 0$,

$L(-m,\chi) = q^m\sum_{a=1}^{q}\chi(a)\zeta(-m,\frac{a}{q})$, we have the following corollary of our Theorem 1.

**Corollary 2**: We have for any integer $m \geq 0$,

$$L(-m,\chi)=\sum_{k=1}^{m}\binom{m}{k}\zeta(-k)q^{k}\left(\sum_{a=1}^{q}\chi(a)a^{m-k}\right)+(\zeta(0)+1)\left(\sum_{a=1}^{q}\chi(a)a^{m}\right)-\tfrac{1}{q(m+1)}\left(\sum_{a=1}^{q}\chi(a)a^{m+1}\right).$$

Thus $L(-m,\chi)$ for $m\geq 0$ is a linear combination of $S(-n,\chi)=\sum_{a=1}^{q}\chi(a)a^{n}$ for $1\leq n\leq m+1$ with rational coefficients.

**Note**: The above statement of Theorem 1 concerning Hurwitz zeta function is the same as the already known result from the theory of Bernoulli polynomials, namely $B_{m}(\alpha)=\sum_{k=0}^{m}\binom{m}{k}B_{k}\alpha^{m-k}$ with m replaced by $m+1$. Here $B_{n}$'s are Bernoulli numbers and $B_{n}(\alpha)$'s are Bernoulli polynomials defined by $\frac{ze^{\alpha z}}{e^{z}-1}=\sum_{n=0}^{\infty}\frac{B_{n}(\alpha)}{n!}z^{n}$ for $|z|<2\pi$ and $B_{n}=B_{n}(1)$ for $n\geq 0$. This is evident in view of the already known formula $\zeta(-n,\alpha)=-\frac{B_{n+1}(\alpha)}{n+1}$ for $n\geq 0$. In fact, it is evident that the two facts, namely,

I) Our Theorem 1 above and Lemma 4 I) below together can be used to develop the theory of Bernoulli polynomials $B_{n}(\alpha)$ (or the theory of Euler polynomials $E_{n}(\alpha)$) by just noting

$B_{n+1}(\alpha)=-(n+1)\zeta(-n,\alpha)$ for any integer $n\geq o$.

**Note**: 1) Putting $\alpha=1$, we get for any integer $m\geq 1$, $\sum_{0\leq k\leq m-1}\binom{m}{k}\zeta(-k)+1=\frac{1}{m+1}$, where empty sum denotes zero. Substituting $m=1,2,3,\ldots\ldots$ successively, we get the values of $\zeta(0),\zeta(-1),\zeta(-2)\ldots$ successively. Obviously, $\zeta(0),\zeta(-1),\zeta(-2),\ldots$ are all rational numbers.

2) Alternatively, noting that $\zeta(s,\tfrac{1}{2})=(2^{s}-1)\zeta(s)$ and putting $\alpha=\tfrac{1}{2}$ in the statement of the Theorem 1, we get for any integer $m\geq 2$

$$(2^{-m}-2)\zeta(-m)=\sum_{k=0}^{m}\binom{m}{k}\zeta(-k)2^{k-m}+2^{-m}-\tfrac{2^{-m-1}}{m+1}.$$

Substituting $m=2,3,4,\ldots\ldots$, we get the values $\zeta(-1),\zeta(-2),\ldots\ldots$ successively.

3) Alternatively, noting that $\zeta_{1}(s,-\tfrac{1}{2})=(2^{s}-1)\zeta(s)$ and putting $\alpha=-\tfrac{1}{2}$ in the statement

$\zeta_{1}(-m,\alpha)=\sum_{\ell=0}^{m}\binom{m}{\ell}\zeta(-\ell)\alpha^{m-\ell}+\alpha^{m}-\tfrac{\alpha^{m+1}}{m+1}$ for $m\geq 0$, we have

: 7 :

$$(2^{-m}-2)\zeta(-m)=\sum_{\ell=1}^{m-1}\binom{m}{\ell}\zeta(-\ell)(-2)^{\ell-m}+\zeta(o)(-2)^{-m}-\frac{(-2)^{m+1}}{m+1}$$ and substituting $m=2,3,\ldots\ldots\ldots$,

we get the values $\zeta(-1),\zeta(-2),\zeta(-3),\ldots\ldots\ldots\ldots\ldots$ successively.

**Proposition 2 :** For any fixed complex number s, if $0<|\alpha|<1$ we have

$$\alpha^{-s}+\sum_{n\geq 1}\frac{(-1)^n}{n!}s(s+1)\ldots\ldots(s+n-1)\zeta(s+n,\alpha)=0$$

**Remark :** The proof follows on considering $\zeta(s,\alpha+1)$ (for a fixed complex number s) as an analytic function of second variable and by developing $\zeta(s,\alpha+1)$ as a Taylor series around $\alpha$ for $0<|\alpha|<1$ Note that $\zeta(s,\alpha+1)-\zeta(s,\alpha)=-\alpha^{-s}$ and

$$\frac{\partial^n}{\partial\alpha^n}\zeta(s,\alpha)=(-1)^n s(s+1)\ldots\ldots(s+n-1)\zeta(s+n,\alpha).$$

As a consequence of Proposition 2, we state Theorem 2 below, which runs parallel to Theorem 1. In author [5], a more general result than the Theorem 2 below has been obtained by Integration approach, in contrast to proof of Theorem 2 below, which has been obtained by Taylor series approach .

**Theorem 2 :** We have, for any integer $m\geq 1$, $\sum_{k=o}^{m-1}\binom{m}{k}\zeta(-k,\alpha)+\alpha^m=\frac{1}{m+1}$, where $\binom{0}{0}=1$.

**Corollary :** If $\chi(\bmod q)$ is a non-principal Dirichlet character, we have for any integer $m\geq 1$,

$$\sum_{k=o}^{m-1}q^{m-k}\binom{m}{k}L(-k,\chi)+\sum_{a=1}^{q}\chi(a)a^m=0.$$

**Remark :** The corollary follows on noting that $L(-k,\chi)=q^k\sum_{a=1}^{q}\chi(a)\zeta(-k,\frac{a}{q})$ for any integer $k\geq 0$.

**Note :** Taking $m=1,2,3,\ldots\ldots$ successively in the statement of Theorem 2, we get the polynomial expressions for $\zeta(0,\alpha),\zeta(-1,\alpha),\zeta(-2,\alpha),\ldots\ldots$ respectively.

**Proposition 3** : We have

a) $\zeta(s,\alpha) = (2^s - 1)\zeta(s) + \sum_{n\geq 1} \frac{(-1)^n}{n!} s(s+1).......(s+n-1)(2^{s+n} - 1)\zeta(s+n)(\alpha - \frac{1}{2})^n$ for $|\alpha - \frac{1}{2}| < \frac{1}{2}$.

b) $\zeta(s,1+\alpha) = (2^s - 1)\zeta(s) + \sum_{n\geq 1} \frac{(-1)^n}{n!} s(s+1).......(s+n-1)(2^{s+n} - 1)\zeta(s+n)(\alpha + \frac{1}{2})^n$ for $|\alpha + \frac{1}{2}| < \frac{1}{2}$.

c) $\zeta(s,\alpha) = \zeta(s) + \sum_{n\geq 1} \frac{(-1)^n}{n!} s(s+1).......(s+n-1)\zeta(s+n)(\alpha - 1)^n$ for $|\alpha - 1| < 1$.

As a consequence of Proposition 3, we get Theorem 3 below.

**Note** : In Proposition 3, a) follows by considering the Taylor series of $\zeta(s,\alpha)$ about $\alpha = \frac{1}{2}$

and on noting $\zeta(s,\frac{1}{2}) = (2^s - 1)\zeta(s)$ and

b) follows by considering the Taylor series of $\zeta_1(s,\alpha)$ about $\alpha = -\frac{1}{2}$ and on noting that

$\zeta_1(s,-\frac{1}{2}) = (2^s - 1)\zeta(s)$ and

c) follows by considering the Taylor series of $\zeta(s,\alpha)$ about $\alpha = 1$.

**Theorem 3** : We have for any integer $m \geq o$ and any for a complex variable $\alpha$,

I) $\zeta(-m,\alpha) = \sum_{n=o}^{m-1} \binom{m}{n} (2^{n-m} - 1)\zeta(n-m)(\alpha - \frac{1}{2})^n - \frac{(\alpha - \frac{1}{2})^{m+1}}{m+1}$.

II) $\zeta(-m,\alpha) = \alpha^m + \sum_{n=o}^{m-1} \binom{m}{n} (2^{n-m} - 1)\zeta(n-m)(\alpha + \frac{1}{2})^n - \frac{(\alpha + \frac{1}{2})^{m+1}}{m+1}$.

III) $\zeta(-m,\alpha) = \sum_{n=o}^{m} \binom{m}{n} \zeta(n-m)(\alpha - 1)^n - \frac{(\alpha - 1)^{m+1}}{m+1}$.

**Note** : I), II), III) follows from a), b), c) of Proposition 3 respectively.

**Theorem 4 :** For real $\lambda$ and for complex $\alpha$ with, $0 < |\alpha| < 1$ we have

I) $\phi(\lambda,\alpha,s) - \alpha^{-s} = \phi(\lambda,s) + \sum_{n\geq 1} \frac{(-\alpha)^n}{n!} s(s+1)...............(s+n-1)\phi(\lambda,s+n)$ so that for any integer $m \geq 0$

and for real non-integral $\lambda$, we have $\phi(\lambda,\alpha,-m) = \alpha^m + \sum_{k=0}^{m} \binom{m}{k} \phi(\lambda,-k)\alpha^{m-k}$.

II) For real $\lambda$ and complex $\alpha$ with $0 < |\alpha| < 1$ we have

: 9 :

$$(e^{-2\pi i\lambda}-1)\phi(\lambda,\alpha,s)=e^{-2\pi i\lambda}\alpha^{-s}+\sum_{n\geq 1}\frac{(-1)^n}{n!}s(s+1)\ldots\ldots(s+n-1)\phi(\lambda,\alpha,s+n)$$

so that we have for any integer $m\geq 0$,

$$(e^{-2\pi i\lambda}-1)\phi(\lambda,\alpha,-m)=e^{-2\pi i\lambda}\alpha^m+\sum_{n=1}^{m}\binom{m}{n}\phi(\lambda,\alpha,n-m)=e^{-2\pi i\lambda}\alpha^m+\sum_{k=0}^{m-1}\binom{m}{k}\phi(\lambda,\alpha,-k).$$

**Note :** I) follows by considering for $\sigma>1$ and $\lambda$ real, $\phi_1(\lambda,\alpha,s)=\sum_{n\geq 1}e^{2\pi i\lambda n}(n+\alpha)^{-s}$ and imitating the corresponding proof for $\zeta(s,\alpha)$ of author [3]. Also note that for non-integral real $\lambda$, $\phi(\lambda,\alpha,-m)$ is a polynomial in $\alpha$ of degree m unlike $\zeta(-m,\alpha)$, which is a polynomial in $\alpha$ of degree m+1. This is so, because $\phi(\lambda,s)$ is analytic at $s=1$, if $\lambda$ is non-integral, real number unlike $\zeta(s)$.

II) follows on considering $\phi(\lambda,\alpha+1,s)$ as analytic function of the second variable and on expanding it about $\alpha$ as a Taylor series for $0<|\alpha|<1$. Note that $\phi(\lambda,\alpha+1,s)-\phi(\lambda,\alpha,s)=(e^{-2\pi i\lambda}-1)\phi(\lambda,\alpha,s)$ and $\frac{\partial^n}{\partial\alpha^n}\phi(\lambda,\alpha,s)=(-1)^n s(s+1)\ldots\ldots(s+n-1)\phi(\lambda,\alpha,s+n)$ for any integer $n\geq 1$.

Next, we give the power series expression of $Z_r(s,\alpha)$ for any integer $r\geq 1$.

We define for $\sigma>r$, $\zeta_r(s)=\sum_{n_1,n_2,\ldots,n_r\geq 0}(n_1+n_2+\ldots+n_r)^{-s}$, with the summation over those r-tuples $(n_1,n_2,\ldots,n_r)$ with $\sum_{\ell=1}^{r}n_i\geq 1$; and its analytic continuation.

**Theorem 5** : We have for any integer $r\geq 1$ and for complex $\alpha$ with, $0<|\alpha|<1$,

$$Z_r(s,\alpha)=\alpha^{-s}+\zeta_r(s)+\sum_{n\geq 1}\frac{(-\alpha)^n}{n!}s(s+1)\ldots\ldots(s+n-1)\zeta_r(s+n)$$

and if $\zeta_r(s)=\frac{1}{(r-1)!}\sum_{j=1}^{r-1}A_j\zeta(s-j)$, where $A_j$'s are constants, we have for any integer $m\geq 0$ and for $r\geq 2$,

$$Z_r(-m,\alpha)=\alpha^m+\sum_{n=0}^{m}\binom{m}{n}\zeta_r(n-m)\alpha^n+\frac{1}{(r-1)!}\sum_{j=1}^{r-1}\frac{(-1)^{j+1}A_j\alpha^{m+j+1}}{(j+1)\binom{m+j+1}{j+1}}.$$

Next, we state our Theorem 6, which evaluates $L(1,\chi)$, where $\chi\pmod q$ is a primitive Dirichlet character.

**Theorem 6**: For any integer $q > 1$, let $\chi \pmod q$ be a primitive Dirichlet character.

I) If $\chi \pmod q$ is even, then $L(1,\chi) = -\frac{1}{\tau(\overline{\chi})} \sum_{a=1}^{q} \overline{\chi}(a) \log \sin \frac{\pi a}{q}$.

II) If $\chi \pmod q$ is odd, then $L(1,\chi) = -\frac{\pi i}{q\tau(\overline{\chi})} \sum_{a=1}^{q} a \overline{\chi}(a)$.

Next, we state a few lemmas, which along with our Theorem 1 or Theorem 2 answer all the questions raised above.

**Lemma 1**: We have the following functional equations.

I) $\zeta(s) = 2^s \pi^{s-1} \Gamma(1-s) \sin \frac{\pi s}{2} \zeta(1-s)$.

II) For any integer $q \geq 1$, if $\chi \pmod q$ is an even, primitive character, then we have

$$L(s,\chi) = q^{-s} 2^s \pi^{s-1} \Gamma(1-s) \tau(\chi) \sin \frac{\pi s}{2} L(1-s, \overline{\chi})$$

III) For any integer $q \geq 1$, if $\chi \pmod q$ is an odd, primitive Dirichlet character, then we have

$$L(s,\chi) = -i q^{-s} 2^s \pi^{s-1} \Gamma(1-s) \tau(\chi) \cos \frac{\pi s}{2} L(1-s, \overline{\chi}).$$

From Lemma 1 follows Lemma 2.

**Lemma 2**: We have

I) $\zeta(2m) = \frac{(-1)^m \pi^{2m} 2^{2m-1}}{(2m-1)!} \zeta(1-2m)$ for $m \geq 1$.

II) We have for $m \geq 1$, if $\chi \pmod q$ is a primitive, even character, then

$$L(2m,\chi) = \frac{(-1)^m 2^{2m-1} q^{1-2m} \pi^{2m}}{(2m-1)! \tau(\overline{\chi})} L(1-2m, \overline{\chi})$$

III) we have for $m \geq 0$, if $\chi \pmod q$ is primitive, odd Dirichlet character, then

$$L(1+2m,\chi) = \frac{i(-1)^m 2^{2m} \pi^{2m+1} q^{-2m}}{(2m)! \tau(\overline{\chi})} L(-2m, \overline{\chi}).$$

**Lemma 3**: I) We have for $m \geq 1$, $\zeta(2m+1) = \frac{(-1)^m 2^{2m+1} \pi^{2m}}{(2m)!} \zeta'(-2m)$

<pre>                                : 11 :</pre>

$$= \frac{(-1)^m 2^{2m} \pi^{2m}}{(2m)!} \left( \zeta'(-2m, \tfrac{1}{2}) + \zeta'_1(-2m, -\tfrac{1}{2}) \right).$$

II) If $\chi(\bmod q)$ is an even, primitive character, then we have for $m \geq 0$,

$$L(2m+1, \chi) = \frac{(-1)^m 2^{2m+1} \pi^{2m} q^{-2m}}{\tau(\overline{\chi})(2m)!} L'(-2m, \overline{\chi}) = \frac{(-1)^m 2^{2m+1} \pi^{2m}}{\tau(\overline{\chi})(2m)!} \sum_{a=1}^{q} \overline{\chi}(a) \zeta'(-2m, \tfrac{a}{q}).$$

III) If $\chi(\bmod q)$ is an odd, primitive character, then we have for $m \geq 1$.

$$L(2m, \chi) = \frac{i(-1)^{m+1} 2^{2m} \pi^{2m-1} q^{1-2m}}{\tau(\overline{\chi})(2m-1)!} L'(1-2m, \overline{\chi}) = \frac{i(-1)^{m+1} 2^{2m} \pi^{2m-1}}{\tau(\overline{\chi})(2m-1)!} \sum_{a=1}^{q} \overline{\chi}(a) \zeta'(1-2m, \tfrac{a}{q}).$$

**Sketch of the proof of Lemma 3 :** As an illustration, we shall sketch II). Let $\chi(\bmod q)$ be an even, primitive character so that we have, by Lemma 1, $L(s, \chi) = q^{-s} 2^s \pi^{s-1} \Gamma(1-s) \tau(\chi) \sin\frac{\pi s}{2}. L(1-s, \overline{\chi})$.

Consider this functional equation for $\sigma < o$. For $s = -2m$, where $m \geq o$ is an integer, in view of the presence of the factor $\sin\frac{\pi s}{2}$, we have $L(-2m, \chi) = 0$.

Thus, we have $\lim_{s \to -2m} \frac{L(s,\chi)}{\sin\frac{\pi s}{2}} = \lim_{s \to -2m} \frac{L'(s,\chi)}{\frac{\pi}{2}\cos\frac{\pi s}{2}} = \frac{L'(-2m,\chi)}{\frac{\pi}{2}\cos \pi m} = (-1)^m \frac{2}{\pi} L'(-2m, \chi)$

by letting $s \to -2m$ through real values of s in view of L'Hospital's rule.

Thus, $\lim_{s \to -2m} \frac{L(s,\chi)}{\sin\frac{\pi s}{2}} = \lim_{s \to -2m} q^{-s} 2^s \pi^{s-1} \Gamma(1-s) \tau(\chi) L(1-s, \overline{\chi})$

That is, $(-1)^m \frac{2}{\pi} L'(-2m, \chi) = q^{2m} 2^{-2m} \pi^{-2m-1} \Gamma(1+2m) \tau(\chi) L(1+2m, \overline{\chi}).$

This prove II).

**Note** : We have seen that for any integer $n \geq o$, $\zeta(-n), \zeta(-n, \alpha)$ and $L(-n, \chi)$ are explicitly computable, as the power series of $\zeta_1(-n, \alpha)$ is a polynomial in $\alpha$. Hence, in the light of Lemma 2, for any integer $m \geq 1$, $\zeta(2m)$ is computable, when $m \geq 1$, and so also $L(n, \chi)$, when $n \geq 1$ and $\chi(\bmod q)$ are of the same parity. However, unlike $\zeta_1(-n, \alpha)$, the power series (in $\alpha$) of $\zeta_1'(-n,\alpha)$ is <u>Not a polynomial</u> and thus cannot be summed up. See Theorem 2 of author [4] or entry 28 (b) of Berndt [1]. Note $\zeta'(-n,\alpha) = -\alpha^n \log \alpha + \zeta_1'(-n,\alpha)$. Thus $\zeta'(-n,\alpha)$ cannot be summed up. In author [5], $\zeta(-n,\alpha)$ and $\zeta'(-n,\alpha)$ for $n \geq 0$ have been obtained by integration approach. There also, the expression for $\zeta'(-n,\alpha)$ contains a term involving an integral, which cannot be evaluated. Consequently, $\zeta'(-n)$, $L'(-n,\chi)$ cannot be evaluated and thus, in the light of Lemma 3, $\zeta(2n+1)$ for $n \geq 1$ cannot computed; so also $L(n,\chi)$, when $n > 1$ and $\chi \pmod{q}$ are of opposite parity, cannot be computed. However, rapidly convergent series expressions for them can be given. See Theorem 4 of author [4].

**Lemma 4**: We have for $\sigma < 0$ and for $0 < \alpha \leq 1$,

I) $\zeta(s,\alpha) = 2^s \pi^{s-1} \Gamma(1-s) \sum_{n \geq 1} \sin(\frac{\pi s}{2} + 2\pi n\alpha) n^{s-1}$

so that for any integer $m \geq 0$, we have $\zeta(-m,\alpha) = 2^{-m} \pi^{-m-1} m! \sum_{n \geq 1} \frac{\sin(\frac{\pi m}{2} + 2\pi n\alpha)}{n^{m+1}} = m! \sum_{|n| \geq 1} \frac{e^{2\pi i n\alpha}}{(2\pi i n)^{m+1}}$

II) We have for $\sigma < 0$ and for $0 < \lambda < 1$ and for $0 < \alpha \leq 1$,

$e^{2\pi i \lambda \alpha} \phi(\lambda,\alpha,s) = i(2\pi)^{s-1} \Gamma(1-s) e^{\frac{-\pi i s}{2}} \sum_{n=-\infty}^{\infty} e^{-2\pi i n\alpha} (n+\lambda)^{s-1}$

so that $\phi(\lambda,\alpha,-m) = i(2\pi)^{-m-1} m! e^{\frac{\pi i m}{2}} \sum_{n=-\infty}^{\infty} \frac{e^{-2\pi i(n+\lambda)\alpha}}{(n+\lambda)^{m+1}}$

**Note :** Lemma 4 I) is the Hurwitz's formula written in suitable form. The same formula can be written as $\zeta(s,\alpha) = \Gamma(1-s) \sum_{|n| \geq 1} e^{2\pi i n\alpha} (2\pi i n)^{s-1}$, where the principal value of logarithm is under consideration. Note that II) gives the Fourier series of $e^{2\pi i \lambda \alpha} \phi(\lambda,\alpha,-m)$ as a function of $\alpha$ on unit interval for a fixed real $\lambda$.

Next, we state our Proposition 4, which in particular, explains why $L(2m+1)$ or $L_1(2m+1)$ for any integer $m \geq 1$ is computable and why $L(2m)$ or $L_1(2m)$ is not computable.

**Proposition 4 :** We have, for any integers $m \geq 1$ and $q \geq 1$,

I) $\zeta(-2m, \frac{1}{q}) = (-1)^m 2^{-2m} \pi^{-2m-1} (2m)! \cdot \sum_{1 \leq r \leq \frac{q-1}{2}} \sin \frac{2\pi r}{q} \left( \sum_{n \equiv r \pmod{q}} - \sum_{n \equiv q-r \pmod{q}} \right) n^{-2m-1}$

II) $\zeta(1-2m, \frac{1}{q}) = (-1)^m 2^{1-2m} \pi^{-2m} (2m-1)! \sum_{r=1}^{q} \cos \frac{2\pi r}{q} \sum_{n \equiv r \pmod{q}} n^{-2m}$

**Proof of Proposition 4:** We have for $\sigma < 0$, $\zeta(s, \alpha) = 2^s \pi^{s-1} \Gamma(1-s) \sum_{n \geq 1} \sin(\frac{\pi s}{2} + 2\pi n \alpha) n^{s-1}$.

For any integers $m \geq 1$ and $q \geq 1$,

$\zeta(-2m, \frac{1}{q}) = 2^{-2m} \pi^{-2m-1} \Gamma(1+2m) \sum_{n \geq 1} \sin(-\pi m + \frac{2\pi n}{q}) n^{-2m-1}$

$= (-1)^m 2^{-2m} \pi^{-2m-1} (2m)! \sum_{n \geq 1} (\sin \frac{2rn}{q}) n^{-2m-1} = (-1)^m 2^{-2m} \pi^{-2m-1} (2m)! \sum_{r=1}^{q} \sin \frac{2\pi r}{q} \sum_{n \equiv r \pmod q} n^{-2m-1}$

$= (-1)^m 2^{-2m} \pi^{-2m-1} (2m)! \sum_{1 \leq r \leq \frac{q}{2}} \sin \frac{2\pi r}{q} \left( \sum_{n \equiv r \pmod q} - \sum_{n \equiv q-r \pmod{}} \right) n^{-2m-1}.$

In Proposition 2 (I), taking $q = 4$, we get

$\zeta(-2m, \frac{1}{4}) = (-1)^m 2^{-2m} \pi^{-2m-1} (2m)! \left( \sum_{n \equiv 1 \pmod 4} - \sum_{n \equiv 3 \pmod 4} \right) n^{-2m-1}.$

Note that $L(2m+1) = \left( \sum_{n \equiv 1 \pmod 4} - \sum_{n \equiv 3 \pmod 4} \right) n^{-2m-1}$ and thus we have the following corollary:

**Corollary :** we have $L(2m+1) = \frac{(-1)^m 2^{2m} \pi^{2m+1}}{(2m)!} \zeta(-2m, \frac{1}{4})$ so that $L(2m+1)$ is a rational multiple of

$\pi^{2m+1}$. Note that evaluation of $\zeta(2n+1) = \sum_{n \geq 1} n^{-2m-1}$ amounts to the evaluation of

$\left( \sum_{n \equiv 1 \pmod 4} + \sum_{n \equiv 3 \pmod 4} \right) n^{-2m-1}$. Instead, we can compute $\left( \sum_{n \equiv 1 \pmod 4} - \sum_{n \equiv 3 \pmod 4} \right) n^{-2m-1}$. Similarly, instead of being able to compute $\left( \sum_{n \equiv 1 \pmod 3} + \sum_{n \equiv 2 \pmod 3} \right) n^{-2m-1}$, we can compute $\left( \sum_{n \equiv 1 \pmod 3} - \sum_{n \equiv 2 \pmod 3} \right) n^{-2m-1}$ as rational multiple of $\zeta(-2m,\frac{1}{3})\pi^{2m+1}$ Similar is the case with $\left(\sum_{n\equiv 1(\bmod 6)} - \sum_{n\equiv 5(\bmod 6)}\right)n^{-2m-1}$, which is a

rational multiple of $\zeta(-2m,\frac{1}{6})\pi^{2m+1}\sin\frac{\pi}{3}$. Similarly, we can compute

$\left(\sum_{n\equiv 1(\bmod 8)} - \sum_{n\equiv 5(\bmod 8)} + \sum_{n\equiv 3(\bmod 8)} - \sum_{n\equiv 7(\bmod 8)}\right)n^{-2m-1}$ and so on. These facts give rise to the following conjucture.

**Conjucture** : For any integer $n \geq 1$, $\zeta(2n+1)$ is an algebraic number multiple of $\pi^{2n+1}$.

Next, we give the proof of our main Theorem 1. We give apparently three different methods, though, in principle, all of them are the same as they are based on the fact that $\frac{\partial}{\partial \alpha}\zeta(s,\alpha) = -s\zeta(s+1,\alpha)$. We also use the fact that $\lim_{s \to 0} s\zeta(s+1,\alpha) = 1$ for any $\alpha$ with $0 < \alpha \leq 1$.
In particular, $\lim_{s \to 0} s\zeta(s+1) = 1$.

**Proof of Theorem 1** :

**Method 1** : We have for any real variable $\alpha$ with $[0,1]$ and for any integer $m \geq 1$.

$\zeta_1(-m,\alpha) = \zeta(-m) + \left(\sum_{1\leq n \leq m} + \sum_{n=m+1} + \sum_{n \geq m+2}\right)\frac{(-\alpha)^n}{n!}s(s+1).......(s+n-1)\zeta(s+n)$, with $s = -m$.

With $s = -m$, consider the single term $\sum_{n=m+1}\frac{(-\alpha)^n}{n!}s(s+1).........(s+n-1)\zeta(s+n)$

$= \frac{(-\alpha)^{m+1}}{(m+1)!}(s(s+1).......(s+m-1))(s+m)\zeta(s+m+1)$

$= \frac{(-\alpha)^{m+1}}{(m+1)!}\left((-1)^m m(m-1).......1\right)\lim_{s \to -m}(s+m)\zeta(s+m+1) = -\frac{\alpha^{m+1}}{(m+1)!}m! = -\frac{\alpha^{m+1}}{m+1}$,

where we have used the fact that $\lim_{w \to 0} w\zeta(w+1) = 1$.

Next with $s = -m$, consider the summation $\sum_{n \geq m+2}\frac{(-\alpha)^n}{n!}s(s+1).......(s+n-1)\zeta(s+n)$.

Note that each term in the summation has zero as a factor.

Thus $\sum_{n \geq m+2} \frac{(-\alpha)^n}{n!} s(s+1).......(s+n-1)\zeta(s+n) = 0$ with $s = -m$.

Thus we have $\zeta_1(-m,\alpha) = \zeta(-m) + \sum_{n=1}^{m} \frac{(-\alpha)^n}{n!} s(s+1).......(s+n-1)\zeta(s+n) - \frac{\alpha^{m+1}}{m+1}$, with $s = -m$.

This gives, $\zeta_1(-m,\alpha) = \zeta(-m) + \sum_{n=1}^{m} \frac{(-\alpha)^n}{n!}(-1)^n m(m-1).......(m-n+1)\zeta(n-m) - \frac{\alpha^{m+1}}{m+1}$

Thus $\zeta_1(-m,\alpha) = \zeta(-m) + \sum_{n=1}^{m} \binom{m}{m-n}\zeta(n-m)\alpha^n - \frac{\alpha^{m+1}}{m+1} = \zeta(-m) + \sum_{k=0}^{m-1} \binom{m}{k}\zeta(-k)\alpha^{m-k} - \frac{\alpha^{m+1}}{m+1}$.

$= \sum_{k=0}^{m} \binom{m}{k}\zeta(-k)\alpha^{m-k} - \frac{\alpha^{m+1}}{m+1}$. Since $\zeta(-m,\alpha) = \alpha^m + \zeta_1(-m,\alpha)$, we get the required result.

**Method 2 :** We have for a fixed integer $m \geq 1$ and for $\alpha$ on the interval $[0,1]$,

$\zeta_1(-m,\alpha) = \zeta(-m) + \sum_{n=1}^{m} \frac{(-\alpha)^n}{n!}s(s+1)........(s+n-1)\zeta(s+n) + \frac{(-\alpha)^{m+1}}{(m+1)!}s(s+1).......(s+m)\zeta(s+m+1,\theta\alpha)$,

where $s = -m$ and $\theta = \theta(m)$ depends on m with $0 < \theta < 1$.

Thus $\zeta(-m,\alpha) = \zeta(-m) + \sum_{n=1}^{m} \frac{(-\alpha)}{n!}(-1)m(m-1).......(m-n+1)\zeta(n-m)$

$+ \frac{(-\alpha)^{m+1}}{(m+1)!}(-1)^m (m(m-1)(m-2)..........1) \lim_{s \to -m} (s+m)\zeta(s+m+1,\theta\alpha)$.

$= \zeta(-m) + \sum_{n=1}^{m} \binom{m}{n}\zeta(n-m)\alpha^n - \frac{\alpha^{m+1}}{m+1}$, where we have used the fact that

$\lim_{s \to 0} s\zeta(s+1,\alpha) = 1$ for any $\alpha$ with $0 < \alpha \leq 1$. Thus, we have

$\zeta_1(-m,\alpha) = \sum_{n=0}^{m} \binom{m}{m-n}\zeta(n-m)\alpha^n - \frac{\alpha^{m+1}}{m+1} = \sum_{k=0}^{m} \binom{m}{k}\zeta(-k)\alpha^{m-k} - \frac{\alpha^{m+1}}{m+1}$.

As $\zeta(-m,\alpha) = \alpha^m + \zeta_1(-m,\alpha)$, we get the required result. This result is apparently valid for $|\alpha| < 1$.

However for any integer $k \geq 1$, if we start with $\zeta_k(s,\alpha)$ instead of $\zeta_1(s,\alpha)$ and proceed, we get the same result for $|\alpha| < k$. However $k \geq 1$ is arbitrary and thus we have the above result for $\zeta(-m,\alpha)$ for any complex number $\alpha$..

**Method 3**: We know that for $\sigma < 0, \zeta(s,\alpha) = \alpha^{-s} + \zeta_1(s,\alpha)$ is a continuous function of $\alpha$ on the interval $[0,1]$, for a fixed complex number s. Also note that for $\sigma < 0, \zeta(s,o) = \zeta_1(s,o) = \zeta(s)$.

We begin with the fact $\zeta(o,\alpha) = \frac{1}{2} - \alpha$ and use the fact that $\frac{\partial}{\partial \alpha}\zeta(s,\alpha) = -s\zeta(s+1,\alpha)$

so that $\zeta(o,\alpha) = \frac{1}{2} - \alpha$ gives

$$\int_o^\alpha \frac{\partial}{\partial \beta}(-\zeta(-1,\beta))d\beta = \int_o^\alpha (\tfrac{1}{2} - \beta)d\beta.$$

Thus $\zeta(-1,\alpha) - \zeta(-1) = \frac{\alpha^2}{2} - \frac{\alpha}{2}$, so that $\zeta(-1,\alpha) = \frac{\alpha^2}{2} - \frac{\alpha}{2} + \zeta(-1)$. Hence

$$\int_o^\alpha \zeta(-1,\beta)d\beta = \int_o^\alpha \left(\frac{\beta^2}{2} - \frac{\beta}{2} + \zeta(-1)\right) d\beta$$

That is, $\zeta(-2,\alpha) - \zeta(-2) = (-2)\left(\frac{\alpha^3}{3!} - \frac{\alpha^2}{4} + \zeta(-1)\alpha\right)$.

We integrate once more with respect to $\alpha$, and we continue this process.

Thus, we get for any integer $m \geq 0$

$$\zeta(-m,\alpha) = \sum_{k=1}^m \binom{m}{k} \zeta(-k)\alpha^{m-k} + (\zeta(0)+1)\alpha^m - \frac{\alpha^{m+1}}{m+1},$$ where empty sum stands for zero.

**Proof of Theorem 2**: From Proposition 2, we have $\alpha^{-s} + \sum_{n\geq 1} \frac{(-1)^n s(s+1)\dots(s+n-1)\zeta(s+n,\alpha)=o}{n!} = 0$

To begin with, let $0 < \alpha \leq 1$. We have for $s = -m$, where $m \geq 0$ is an integer,

$$\alpha^m + \left(\sum_{n=1}^m + \sum_{n=m+1} + \sum_{n\geq m+2}\right) \frac{(-1)^n}{n!} s(s+1)\dots(s+n-1)\zeta(s+n,\alpha) = 0$$

As in the proof of Theorem 1, with $s = -m$, the term with $n = m+1$ gives $-\frac{1}{m+1}$ on using the fact

$\lim_{s \to 0} s\zeta(s+1,\alpha) = 1$ for $0 < \alpha \leq 1$ and the summation $\sum_{n\geq m+2} = 0$, as each term in the summation

has zero as a factor. Thus we are left with $\alpha^m + \sum_{n=1}^m \frac{(-1)^n}{n!}(-m)(-m+1)\dots(-m+n-1)\zeta(n-m,\alpha) - \frac{1}{m+1} = 0$

That is $\alpha^m + \sum_{n=1}^m \frac{(-1)^n}{n!}\binom{m}{m-n}\zeta(n-m,\alpha) - \frac{1}{m+1} = 0$.

That is $\alpha^m + \sum_{k=0}^{m-1} \binom{m}{k} \zeta(-k,\alpha) - \frac{1}{m+1} = 0$, on putting $k = m - n$.

This completes the proof of Theorem 2.

**Sketch of proof of Theorem 5**: We have for $\sigma > r$, $Z_r(s,\alpha) - \alpha^{-s} = \sum_{N_r \geq 1}(n_1 + n_2 \ldots + n_r + \alpha)^{-s}$, where $N_r = \sum_{i=1}^{r} n_i$.

Thus $Z_r(s,\alpha) - \alpha^{-s} = \sum_{N_r \geq 1}(N_r + \alpha)^{-s} = \sum_{N_r \geq 1} N_r^{-s}(1 + \frac{\alpha}{N_r})^{-s}$

$= \sum_{N_r \geq 1} N_r^{-s} - \alpha s \sum_{N_r \geq 1} N_r^{-s-1} + \frac{\alpha^2}{2!} s(s+1) \sum_{N_r \geq 1} N_r^{-s-2} \ldots$

$= \zeta_r(s) + \sum_{n \geq 1} \frac{(-\alpha)^n}{n!} s(s+1)\ldots(s+n-1)\zeta_r(s+n)$ and continue on similar lines as in author [3] as in the case of $\zeta(s,\alpha)$. This gives for $0 < |\alpha| < 1$,

$Z_r(s,\alpha) = \alpha^{-s} + \zeta_r(s) + \sum_{n \geq 1} \frac{(-\alpha)^n}{n!} s(s+1)\ldots(s+n-1)\zeta_r(s+n)$

Next, we evaluate $Z_r(-m,\alpha)$ for any integer $m \geq 0$ and for $r \geq 2$. We have $Z_r(-m,\alpha)$

$= \alpha^m + \zeta_r(-m) + \left(\sum_{n=1}^{m} + \sum_{n=m+1} + \sum_{n=m+2}^{m+r} + \sum_{n \geq m+r+1}\right) \frac{(-\alpha)^n}{n!} s(s+1)\ldots(s+n-1)\zeta_r(s+n)$ with $s = -m$.

Consider $\sum_{n \geq m+r+1} \frac{(-\alpha)^n}{n!} s(s+1)\ldots(s+n-1)\zeta_r(s+n)$ with $s = -m$.

This is $= \sum_{n \geq m+r+1} \frac{(-\alpha)^n}{n!}(-m+1)\ldots(-m+n-1)\zeta_r(n-m)$. Note that as $n \geq m+r-1$, each term in the summation has zero as a factor and $\zeta_r(n-m) = \frac{1}{(r-1)!}\sum_{j=1}^{r-1} A_j \zeta(n-m-j)$ is finite for each $n \geq m+r+1$.

Thus summation $\sum_{n \geq m+r+1} \frac{(-\alpha)^n}{n!} s(s+1)\ldots(s+n-1)\zeta_r(s+n) = 0$ with $s = -m$.

Next consider $\sum_{n=m+1} \frac{(-\alpha)^n}{n!} s(s+1)\ldots(s+n-1)\zeta_r(s+n)$ with $s = -m$.

This is $= \frac{(-\alpha)^{m+1}}{(m+1)!}(-m)(-m+1)\ldots\ldots\ldots 1.0.\zeta_r(1)$. Note that $\zeta_r(1) = \frac{1}{(r-1)}\sum_{j=1}^{r-1} A_j \zeta(1-j)$ is finite.

Thus $\sum_{n=m+1} = 0$.

Consider $\sum_{n=1}^{m} \frac{(-\alpha)^n}{n!} s(s+1)\ldots\ldots\ldots(s+n-1)\zeta_r(s+n)$ with $s = -m$.

This summation $= \sum_{n=1}^{m} \frac{(-\alpha)^n}{n!}(-m)(-m+1)\ldots\ldots\ldots(-m+n-1)\zeta_r(n-m) = \sum_{n=1}^{m} \binom{m}{n} \zeta_r(n-m)\alpha^n$.

Next consider for $s = -m$, $\sum_{n=m+2}^{m+r} \frac{(-\alpha)^n}{n!} s(s+1)\ldots\ldots\ldots(s+n-1)\zeta_r(s+n)$

$= \sum_{n=m+2}^{m+r} \frac{(-\alpha)^n}{n!}(-m)(-m+1)\ldots\ldots\ldots(n-m-1)\zeta_r(n-m) = \sum_{n=m+2}^{m+r} \frac{(-\alpha)^n}{n!}(-m)(-m+1)\ldots\ldots\ldots(n-m-1)\zeta_r(n-m)$

$= \sum_{n=m+2}^{m+r} \frac{(-\alpha)^n}{n!}(-m)(-m+1)(-2)(-1)0.1.2\ldots\ldots\ldots(n-m-1)\zeta_r(n-m) = \sum_{n=m+2}^{m+r} \frac{\alpha^n(-1)^{m+n}}{n!} m!(n-m-1)!(0.\zeta_r(n-m))$

$= \sum_{n=m+2}^{m+r} \frac{(-1)^{m+n}}{(n-m)} \alpha^n \left(\frac{m!(n-m)!}{n!}\right)\left(0.\frac{1}{(r-1)!}\sum_{j=1}^{r-1} A_j \zeta(n-m-j)\right) = \frac{1}{(r-1)!} \sum_{n=m+2}^{m+r} \frac{(-1)^{m+n}\alpha^n}{(n-m)\binom{n}{m}}\left(\sum_{j=1}^{r-1} A_j(0.\zeta(n-m-j))\right)$

$= \frac{1}{(r-1)!} \sum_{j=1}^{r-1} \sum_{n=m+j+1} \frac{(-1)^{m+n} A_j \alpha^n}{(n-m)\binom{n}{m}} = \frac{1}{(r-1)!}\sum_{j=1}^{r-1} \frac{(-1)^{2m+j+1} A_j \alpha^{m+j+1}}{(j+1)\binom{m+j+1}{m}} = \frac{1}{(r-1)!}\sum_{j=1}^{r-1} \frac{(-1)^{j+1} A_j \alpha^{m+j+1}}{(j+1)\binom{m+j+1}{m}}$

Here, we have used the fact that $0.\zeta(1) = \lim_{s \to 0} s\zeta(s+1) = 1$.

Thus $Z_r(-m,\alpha) = \alpha^m + \sum_{n=0}^{m} \binom{m}{n} \zeta_r(n-m)\alpha^n + \frac{1}{(r-1)!}\sum_{j=1}^{r-1} \frac{(-1)^{j+1} A_j}{j+1} \frac{\alpha^{m+j+1}}{\binom{m+j+1}{m}}$.

**Proof of Theorem 6:** If $\chi(\mod q)$ is an odd primitive character, then using Lemma 2, we have

$$L(1,\chi) = \frac{\pi i}{\tau(\chi)} L(0,\overline{\chi}) = \frac{\pi i}{\tau(\chi)} \sum_{a=1}^{q} \overline{\chi}(a)\zeta(0,\tfrac{a}{q}) = \frac{\pi i}{\tau(\chi)} \sum_{a=1}^{q} \overline{\chi}(a)(\tfrac{1}{2}-\tfrac{a}{q}) = \frac{-\pi i}{q.\tau(\chi)} \sum_{a=1}^{q} a\overline{\chi}(a).$$

This proves Theorem 6 II).

Next, let $\chi(\mod q)$ be an even primitive character.

Using Lemma 3, we have $L(1,\chi) = \frac{2}{\tau(\chi)} L'(0,\overline{\chi})$. From $L(s,\chi) = q^{-s}\sum_{a=1}^{q} \chi(a)\zeta(s,\tfrac{a}{q})$ on differentiation with respect to s and, on noting that $L(0,\chi) = 0$, we find that $L'(0,\chi) = \sum_{a=1}^{q} \chi(a)\zeta'(0,\tfrac{a}{q})$.

: 19 :

However, we know for $0 < \alpha \leq 1$, $\zeta'(o,\alpha) = \log\frac{\Gamma(\alpha)}{\sqrt{2\pi}}$. Thus we find that

$$L'(o,\bar{\chi}) = \sum_{a=1}^{q} \bar{\chi}(a) \log\frac{\Gamma(\frac{a}{q})}{\sqrt{2\pi}} = \sum_{a=1}^{q} \bar{\chi}(a) \log\Gamma(\frac{a}{q})$$

$$= \sum_{1 \leq a \leq \frac{q-1}{2}} \bar{\chi}(a)\left[\log\Gamma(\frac{a}{q}) + \log\Gamma(1-\frac{a}{q})\right] = \sum_{1 \leq a \leq \frac{q-1}{2}} \bar{\chi}(a) \log\left(\Gamma(\frac{a}{q})\Gamma(1-\frac{a}{q})\right) = \sum_{1 \leq a \leq \frac{q-1}{2}} \bar{\chi}(a) \log\frac{\pi}{\sin\frac{\pi a}{q}}$$

$$= -\sum_{1 \leq a \leq \frac{q-1}{2}} \bar{\chi}(a) \log\sin\frac{\pi a}{q} = -\frac{1}{2}\sum_{a=1}^{q} \bar{\chi}(a) \log\sin\frac{\pi a}{q}.$$

Thus $L(1,\chi) = -\frac{1}{\tau(\chi)} \sum_{a=1}^{q} \bar{\chi}(a) \log\sin\frac{\pi a}{q}$

This completes the proof of Theorem 6.